\documentclass[a4paper]{amsart}
\usepackage{amssymb}
\usepackage{verbatim}
\theoremstyle{plain}
  \newtheorem{thm}{Theorem}[section]
  \newtheorem{lem}[thm]{Lemma}
  \newtheorem{cor}[thm]{Corollary}

  \newtheorem*{obs*}{Observation}
\theoremstyle{definition}

\theoremstyle{remark}
  \newtheorem{rem}[thm]{Remark}
  
  \newtheorem*{ack}{Acknowledgments}

\newcommand{\C}{\mathbb{C}}
\newcommand{\R}{\mathbb{R}}

\newcommand{\Vol}{\operatorname{Vol}}
\newcommand{\arccosh}{\operatorname{arccosh}}
\renewcommand{\Re}{\operatorname{Re}}
\renewcommand{\Im}{\operatorname{Im}}
\numberwithin{equation}{section}
\allowdisplaybreaks
\begin{document}
\title[The colored Jones polynomials and the Alexander polynomial]
{The colored Jones polynomials and the Alexander polynomial of the figure-eight knot}
\author{Hitoshi Murakami}
\address{
Department of Mathematics,
Tokyo Institute of Technology,
Oh-okayama, Meguro, Tokyo 152-8551, Japan
}
\email{starshea@tky3.3web.ne.jp}
\date{\today}
\begin{abstract}
The volume conjecture and its generalization state that the series of certain evaluations of the colored Jones polynomials of a knot would grow exponentially and its growth rate would be related to the volume of a three-manifold obtained by Dehn surgery along the knot.
In this paper, we show that for the figure-eight knot the series converges in some cases and the limit equals the inverse of its Alexander polynomial.
\end{abstract}
\thanks{This research is partially supported by Grant-in-Aid for Scientific
Research (B) (15340019).}
\maketitle
\section{Introduction}
Let $K$ be a knot and $J_N(K;t)$ its colored Jones polynomial corresponding to the $N$-dimensional irreducible representation of $sl_2(\C)$ normalized so that $J_N(U;t)=1$ for the unknot $U$.
The volume conjecture \cite{Murakami/Murakami:ACTAM101} states that
\begin{equation*}
  \lim_{N\to\infty}
  \frac{\log\left|J_N\left(K;\exp(2\pi\sqrt{-1}/N)\right)\right|}{N}
  =
  \frac{v_3}{2\pi}\Vol(S^3\setminus{K}),
\end{equation*}
where $v_3$ is the hyperbolic volume of the ideal regular hyperbolic tetrahedron, and $\Vol$ denotes the simplicial volume.
Note that this conjecture was first proposed by R.~Kashaev \cite{Kashaev:LETMP97} in a different way.
It is generalized by S.~Gukov \cite{Gukov:03} to a relation of the limit
\begin{equation}\label{eq:limit}
  \lim_{N\to\infty}
  \frac{\log J_N\left(K;\exp(a/N)\right)}{N}
\end{equation}
with a fixed complex number $a$ to the $A$-polynomial of $K$
\cite{Cooper/Culler/Gillet/Long/Shalen:INVEM1994}, and the volume and the Chern--Simons invariant of a three-manifold obtained by Dehn surgery along $K$.
See also \cite{Murakami:KYUMJ2004,Murakami/Yokota:2004} about the generalized volume conjecture for the figure-eight knot.
\par
On the other hand, the author proved also in \cite{Murakami:KYUMJ2004} that the  limit \eqref{eq:limit} vanishes for the figure-eight knot if $a$ is real and $|a|<\arccosh(3/2)$ or $a$ is purely imaginary and $|a|<\pi/3$.
S.~Garoufalidis and T.~Le proved \cite[Theorem 2]{Garoufalidis/Le:2005} that for any knot $K$, \eqref{eq:limit} vanishes if $a$ is purely imaginary and sufficiently small.
This shows that the series $\{J_N(K;\exp(a/N))\}_{N=2,3,\dots}$ grows polynomially when $a$ is small.
One may ask whether the series diverges or not.
\par
In this paper we study the `genuine' limit $\lim_{N\to\infty}J_N(E;\exp(a/N))$ for the figure-eight knot $E$ when $a$ is a small complex number, and show that the limit does exist and equals the inverse of its Alexander polynomial.
More precisely we will show the following equality.
\begin{thm}\label{thm}
Let $E$ be the figure-eight knot.
If $a$ is a complex number with $|2\cosh{a}-2|<1$ and $|\Im{a}|<\pi/3$,
then the series $\{J_{N}(E;\exp(a/N))\}_{N=2,3,\dots}$ converges and
\begin{equation*}
  \lim_{N\to\infty}J_{N}\left(E;\exp\frac{a}{N}\right)
  =
  \frac{1}{\Delta(E;\exp{a})},
\end{equation*}
where $\Delta(E;t)=-t+3-t^{-1}$ is the Alexander polynomial of E.
\end{thm}
\begin{rem}
The range $\{a\in\C\mid|2\cosh{a}-2|<1,|\Im{a}|<\pi/3\}$ looks like an `oval' (not a mathematical one) around the origin whose boundary goes through the four points
$\log\left((3+\sqrt{5})/2\right)$, $\pi\sqrt{-1}/3$, $-\log\left((3+\sqrt{5})/2\right)$, and $-\pi\sqrt{-1}/3$ on the Gaussian plane (see Lemma~\ref{lem:cosh}).
The author does not know whether this `oval' is the `circle of convergence' or not.
\end{rem}
\begin{rem}
Note that the inequality $|2\cosh{a}-2|<1$ is equal to $|\Delta(E;\exp{a})-1|<1$.
This may suggest another relation between the colored Jones polynomials and the Alexander polynomial.
\end{rem}
\begin{ack}
This work began when the author was visiting Universit{\'e} de Montpellier II in November, 2004.
The author would like to thank V.~Vershinin for his kind invitation.
\par
Thanks are also due to F.~Nagasato and Y.~Nomura for helpful discussions.
\end{ack}
\section{Proof}
We first recall the formula of the figure-eight knot due to K.~Habiro and Le \cite{Habiro:SURIK00} (see also \cite{Masbaum:AGT03}).
\begin{equation*}
  J_{N}(E;t)
  =
  \sum_{k=0}^{N-1}
  \prod_{j=1}^{k}
  \left(t^{(N+j)/2}-t^{-(N+j)/2}\right)
  \left(t^{(N-j)/2}-t^{-(N-j)/2}\right).
\end{equation*}
If we replace $t$ with $\exp(a/N)$ we have
\begin{equation*}
  J_{N}\left(E;\exp\frac{a}{N}\right)
  =
  \sum_{k=0}^{N-1}f_{N,a}(k)
\end{equation*}
with
\begin{equation*}
  f_{N,a}(k):=\prod_{j=1}^{k}g_{N,a}(j),
\end{equation*}
where
\begin{equation*}
\begin{split}
  g_{N,a}(j)
  &:=
  4\sinh\left(\frac{a(N+j)}{2N}\right)\sinh\left(\frac{a(N-j)}{2N}\right)
  \\
  &=2\cosh{a}-2\cosh\frac{aj}{N}.
\end{split}
\end{equation*}
\par
We first show that $J_{N}\left(E;\exp\frac{a}{N}\right)$ converges.
\begin{lem}
For any complex number $a$ with $|2\cosh{a}-2|<1$ and $|\Im{a}|<\pi/3$,
the series $\{J_{N}\left(E;\exp\frac{a}{N}\right)\}_{N=2,3,\dots}$ converges.
\end{lem}
\begin{proof}
From Lemmas~\ref{lem:cosh_ua_va} and \ref{lem:cosh_ua}, we have the following inequalities for $0<M<N$:
\begin{align*}
  \left|g_{N,a}(j)\right|<&\delta<1
  &\quad&\text{if $0<j<N$,}
  \\
  \left|\frac{g_{M,a}(j)}{g_{N,a}(j)}\right|&<1
  &\quad&\text{if $0<j<M$,}
  \\
  \left|\frac{g_{M,a}(j)}{g_{N,a}(j)}\right|
  &>
  1-\frac{j}{M}
  &\quad&\text{if $0<j<\varepsilon{M}$ for some $\varepsilon>0$,}
\end{align*}
where we put $\delta:=|2\cosh{a}-2|<1$.
So we have
\begin{equation*}
  1>{\delta}^k>
  \left|\frac{f_{M,a}(k)}{f_{N,a}(k)}\right|
  >
  \prod_{j=1}^{\lfloor\varepsilon{M}\rfloor-1}\left(1-\frac{j}{M}\right)
  \prod_{j=\lfloor\varepsilon{M}\rfloor}^{k}
  \left|\frac{f_{M,a}(k)}{f_{N,a}(k)}\right|
\end{equation*}
for $0<k<M<N$, where $\lfloor{x}\rfloor$ is the greatest integer that does not exceed $x$.
\par
Putting $M':=\lfloor\varepsilon{M}\rfloor$ we have
\begin{equation*}
\begin{split}
  &
  \left|
    J_N\left(E;\exp\frac{a}{N}\right)-J_{M}\left(E;\exp\frac{a}{M}\right)
  \right|
  \\
  =&
  \left|
    \sum_{k=0}^{N-1}f_{N,a}(k)
    -
    \sum_{k=0}^{M-1}f_{M,a}(k)
  \right|
  \\
  \le&
  \sum_{k=0}^{M-1}
  \left|
    f_{N,a}(k)-f_{M,a}(k)
  \right|
  +
  \sum_{k=M}^{N-1}\left|f_{N,a}(k)\right|
  \\
  =&
  \sum_{k=0}^{M-1}
  \left|f_{N,a}(k)\right|
  \left(1-\left|\frac{f_{M,a}(k)}{f_{N,a}(k)}\right|\right)
  +
  \sum_{k=M}^{N-1}\left|f_{N,a}(k)\right|
  \\
  <&
  \sum_{k=0}^{M-1}
  {\delta}^k
  \left(
    1-
    \prod_{j=1}^{M'-1}\left(1-\frac{j}{M}\right)
    \prod_{j=M'}^{k}
    \left|\frac{f_{M,a}(k)}{f_{N,a}(k)}\right|
  \right)
  +
  \sum_{k=M}^{N-1}{\delta}^k
  \\
  =&
  \frac{1-{\delta}^N}{1-{\delta}}
  -
  \sum_{k=0}^{M'-1}
  {\delta}^k\prod_{j=1}^{k}\left(1-\frac{j}{M}\right)
  -
  \sum_{k=M'}^{M-1}
  {\delta}^k\prod_{j=1}^{M'-1}\left(1-\frac{j}{M}\right)
  \prod_{j=M'}^{k}
  \left|\frac{f_{M,a}(k)}{f_{N,a}(k)}\right|.
\end{split}
\end{equation*}
From Lemma~\ref{lem:ExpIntegralE}, this is equal to
\begin{multline}\label{eq:J_N-J_M}
  \frac{1-{\delta}^N}{1-{\delta}}
  -
  \frac{M'}{\delta}
  e^{\frac{M'}{\delta}}
  \int_{1}^{\infty}
  e^{-\frac{M'}{\delta}t}
  t^{M'-1}dt
  \\
  -
  \prod_{j=1}^{M'-1}\left(1-\frac{j}{M}\right)
  \sum_{k=M'}^{M-1}
  {\delta}^k\prod_{j=M'}^{k}
  \left|\frac{f_{M,a}(k)}{f_{N,a}(k)}\right|.
\end{multline}
Note that since
\begin{equation*}
  \prod_{j=1}^{M'-1}\left(1-\frac{j}{M}\right)
  \sum_{k=M'}^{M-1}
  {\delta}^k\prod_{j=M'}^{k}
  \left|\frac{f_{M,a}(k)}{f_{N,a}(k)}\right|
  <
  \sum_{k=M'}^{M-1}{\delta}^k
  =
  {\delta}^{M'}
  \frac{1-{\delta}^{M-M'}}{1-{\delta}},
\end{equation*}
the last term in \eqref{eq:J_N-J_M} can be arbitrarily small.
\par
Since
\begin{equation*}
  \int_{1}^{\infty}e^{-\frac{M'}{\delta}t}t^{M'-1}dt
  =
  \int_{1}^{\infty}e^{\left(\log{t}-\frac{t}{\delta}\right)M'}t^{-1}dt,
\end{equation*}
we can apply Laplace's method to study the asymptotic behavior for large $M$:
\begin{equation*}
  \int_{1}^{\infty}e^{\left(\log{t}-\frac{t}{\delta}\right)M'}t^{-1}dt
  \underset{M\to\infty}{\sim}
  \frac{1}{M'}\frac{1}{\frac{1}{\delta}-1}e^{-\frac{M'}{\delta}}
  =
  \frac{\delta}{M'}e^{-\frac{M'}{\delta}}\frac{1}{1-\delta}.
\end{equation*}
(See for example \cite[Chapter~3, \S~7.1]{Olver:1974}.)
Therefore
$\left|
  J_N\left(E;\exp\frac{a}{N}\right)-J_{M}\left(E;\exp\frac{a}{M}\right)
\right|$
can be arbitrarily small if $M$ is sufficiently large, which means that the sequence $\{J_N\left(E;\exp\frac{a}{N}\right)\}_{N=2,3,\dots}$ is a Cauchy sequence and so it converges.
\end{proof}
Now that we know the convergence, we use an {\em inhomogeneous} recursion formula of $J_{N}(E;t)$ to find the limit.
It is known that $J_{N}(E;t)$ satisfies the following formula
\cite[\S 6.2]{Garoufalidis/Le:2003} (see also \cite{Gelca/Sain:JKNOT2004} for a {\em homogeneous} recursion formula).
\begin{equation}\label{eq:recursion}
\begin{split}
  &\phantom{=}
  J_{N}(E;t)
  \\
  &=
  \frac{t^{-N-1}\left(t^N+t\right)\left(t^{2N}-t\right)}
       {t^N-1}
  \\
  &+
  \frac{t^{-2N-2}\left(t^{N-1}-1\right)^2\left(t^{N-1}+1\right)
        \left(t^4+t^{4N}-t^{N+3}-t^{2N+1}-t^{2N+3}-t^{3N+1}\right)}
       {\left(t^N-1\right)\left(t^{2N-3}-1\right)}
  \\
  &\quad\times
  J_{N-1}(E;t)
  \\
  &-
  \frac{\left(t^{N-2}-1\right)\left(t^{2N-1}-1\right)}
       {\left(t^N-1\right)\left(t^{2N-3}-1\right)}
  J_{N-2}(E;t).
\end{split}
\end{equation}
\par
We want to show that the series $\{J_{N-1}(E;\exp\frac{a}{N})\}$ and
$\{J_{N-2}(E;\exp\frac{a}{N})\}$ also converge and both limits coincide with
that of $J_{N}(E;\exp\frac{a}{N})$.
\par
For $l=1$ or $2$, put
\begin{align*}
  g'_{N}(j;l)&:=2\cosh{a\left(1-\frac{l}{N}\right)}-2\cosh{\frac{aj}{N}}
  \\
  \intertext{and}
  f'_{N}(k;l)&:=\prod_{j=1}^{k}{g'}_{N}(j;l)
\end{align*}
so that $J_{N-l}\left(E;\exp{\frac{a}{N}}\right)=\sum_{k=0}^{N-l-1}f'_{N}(k;l)$.
\begin{lem}
The series $\{J_{N-l}(E;\exp\frac{a}{N})\}_{N=2,3,\dots}$ converges and shares the limit with $\{J_{N}(E;\exp\frac{a}{N})\}_{N=2,3,\dots}$.
\end{lem}
\begin{proof}
From Lemma~\ref{lem:cosh_ua_va} and Corollary~\ref{cor}, we have
\begin{equation*}
\begin{split}
  |g'_N(j;l)|
  &=
  2\left|\cosh{a\left(1-\frac{l}{N}\right)}-\cosh\frac{aj}{N}\right|
  \\
  &<
  2\left|\cosh{a\left(1-\frac{l}{N}\right)}-1\right|
  \\
  &<
  2|\cosh{a}-1|
  \\
  &=
  \delta.
\end{split}
\end{equation*}
\par
From Lemma~\ref{lem:cosh_x} there exists a positive number $\varepsilon'$ such that if $j/N<\varepsilon'$ then
\begin{equation*}
\begin{split}
  1
  >
  \left|
    \frac{\cosh{a\left(1-\frac{l}{N}\right)}-\cosh\frac{aj}{N}}
         {\cosh{a}-\cosh\frac{aj}{N}}
  \right|
  >
  1
  -
  \left|
    \frac{a\sinh{a}}{\cosh{a}-1}
  \right|
  \frac{1}{N}.
\end{split}
\end{equation*}
Putting $c:=\left|\frac{a\sinh{a}}{\cosh{a}-1}\right|>0$, we have
\begin{equation*}
  1>\left|\frac{g'_{N}(j;l)}{g_{N}(j)}\right|>1-\frac{c}{N},
\end{equation*}
if $j/N<\varepsilon'$ and so
\begin{equation*}
  1>\left|\frac{f'_{N}(k;l)}{f_{N}(k)}\right|>\left(1-\frac{c}{N}\right)^k
\end{equation*}
if $k/N<\varepsilon'$.
\par
Therefore we have
\begin{equation*}
\begin{split}
  &
  \left|
    J_{N}\left(E;\exp\frac{a}{N}\right)
    -
    J_{N-l}\left(E;\exp\frac{a}{N}\right)
  \right|
  \\
  =&
  \left|
    \sum_{k=0}^{\lfloor{\varepsilon'N}\rfloor-1}
    \left\{f_{N}(k)-f'_{N}(k;l)\right\}
    +
    \sum_{k=\lfloor{\varepsilon'N}\rfloor}^{N-1}
    f_{N}(k)
    -
    \sum_{k=\lfloor{\varepsilon'N}\rfloor}^{N-l-1}
    f'_{N}(k;l)
  \right|
  \\
  <&
  \sum_{k=0}^{\lfloor{\varepsilon'N}\rfloor-1}
  \left|f_N(k)\right|
  \left\{1-\left(1-\frac{c}{N}\right)^k\right\}
  +
  \sum_{k=\lfloor{\varepsilon'N}\rfloor}^{N-1}
  \left|f_{N}(k)\right|
  +
  \sum_{k=\lfloor{\varepsilon'N}\rfloor}^{N-1}
  \left|f'_{N}(k;l)\right|
  \\
  <&
  \sum_{k=0}^{\lfloor{\varepsilon'N}\rfloor-1}
  {\delta}^k\left\{1-\left(1-\frac{c}{N}\right)^k\right\}
  +
  2\sum_{k=\lfloor{\varepsilon'N}\rfloor}^{N-1}
  {\delta}^k
  \\
  =&
  \frac{1-{\delta}^{\lfloor{\varepsilon'N\rfloor}}}{1-\delta}
  -
  \frac{1-{\delta}^{\lfloor{\varepsilon'N\rfloor}}
          \left(1-\frac{c}{N}\right)^{\lfloor{\varepsilon'N\rfloor}}}
       {1-\delta\left(1-\frac{c}{N}\right)}
  +
  2{\delta}^{\lfloor{\varepsilon'N}\rfloor}
  \frac{1-{\delta}^{N-\lfloor{\varepsilon'N}\rfloor}}{1-\delta},
\end{split}
\end{equation*}
which can be arbitrarily small when $N$ is sufficiently large, since $0<\delta<1$.
So the series $\{J_{N-l}\left(E;\exp\frac{a}{N}\right)\}$ ($l=1$ or $2$) converges and its limit is equal to that of $\{J_{N}\left(E;\exp\frac{a}{N}\right)\}$.
\end{proof}
Therefore putting $J_{a}:=\lim_{N\to\infty}J_{N}\left(E;\exp\frac{a}{N}\right)$ and $w:=\exp{a}$, we have from \eqref{eq:recursion}
\begin{equation*}
\begin{split}
  J_{a}
  &=
  \frac{w^{-1}(w+1)\left(w^{2}-1\right)}{w-1}
  \\
  &+
  \frac{w^{-2}(w-1)^2(w+1)\left(1+w^{4}-w-w^{2}-w^{2}-w^{3}\right)}
       {(w-1)\left(w^{2}-1\right)}
  J_{a}
  \\
  &-
  \frac{(w-1)\left(w^{2}-1\right)}
       {(w-1)\left(w^{2}-1\right)}
  J_{a}.
\end{split}
\end{equation*}
So we finally have
\begin{equation*}
  J_{a}=\frac{1}{-w+3-w^{-1}}
\end{equation*}
winch is equal to $1/\Delta(E;\exp{a})$.
\par
This completes the proof of Theorem~\ref{thm}.
\begin{rem}
We used an {\em inhomogeneous} recursion formula for the colored Jones polynomial of the figure-eight knot.
Note that Garoufalidis and Le proved that there always exists a {\em homogeneous} formula for any knot \cite{Garoufalidis/Le:2003}.
\par
The relation between the $A$-polynomial and the Alexander polynomial
\cite[\S 6.3 Proposition]{Cooper/Culler/Gillet/Long/Shalen:INVEM1994}, the AJ-conjecture proposed by Garoufalidis \cite{Garoufalidis:GEOTO2004}, and Theorem~\ref{thm} suggest that for any knot $K$ if the series $\{J_N(K;\exp{a/N})\}_{N=2,3,\dots}$ converges for some $a$, then the limit would be $1/\Delta(K;\exp{a})$ with $\Delta(K;t)$ the Alexander polynomial of $K$.
\par
In \cite{Murakami:INTJM62004} the author proved that for any torus knot $T$, $\lim_{N\to\infty}J_N(T;\exp{a/N})=1/\Delta(T;\exp{a})$ if $a$ is near $2\pi\sqrt{-1}$ and $\Re{a}>0$.
\end{rem}
\begin{rem}
P.~Melvin and H.~Morton \cite{Melvin/Morton:COMMP95} observed the following formal power series:

\begin{equation}\label{eq:MMR}
  J_N(K;\exp{h})=\sum_{j,k\ge0}b_{jk}(K)h^jN^k,
\end{equation}
and conjectured the following (Melvin--Morton--Rozansky conjecture).
\begin{enumerate}
  \item[(i)]
  $b_{jk}(K)=0$ if $k>j$, and
  \item[(ii)]
  $\sum_{j\ge0}b_{jj}(K)(hN)^j=\dfrac{1}{\Delta(K;\exp{hN})}$.
\end{enumerate}
This conjecture was `proved' by L.~Rozansky \cite{Rozansky:COMMP96} non-rigorously, and proved by D.~Bar-Natan and Garoufalidis (\cite{BarNatan/Garoufalidis:INVEM96}).
\par
Replacing $h$ with $a/N$, we have from (i) and (ii)
\begin{align*}
  J_N\left(K;\exp\frac{a}{N}\right)
  &=
  \sum_{j\ge k\ge0}b_{jk}(K)a^jN^{k-j},
  \\
  \intertext{and}
  \sum_{j\ge0}b_{jj}(K)a^j
  &=
  \frac{1}{\Delta(K;\exp{a})}.
\end{align*}
So we may regard Theorem~\ref{thm} as an analytic version of the Melvin--Morton--Rozansky conjecture.
\end{rem}
\section{Appendix}
In this appendix we give several technical lemmas used in the paper.
\begin{lem}\label{lem:cosh}
For a complex number $a=x+y\sqrt{-1}$ with $x$, $y\in\R$, the condition $|2\cosh{a}-2|<1$ is equivalent to the condition $\cosh{x}-\cos{y}<1/2$.
\end{lem}
\begin{proof}
Since
\begin{equation*}
\begin{split}
  |\cosh{a}-1|^2
  &=
  (\cos{y}\cosh{x}-1)^2+\sin^2{y}\sinh^2{x}
  \\
  &=
  (\cos{y}\cosh{x}-1)^2+(1-\cos^2{y})(\cosh^2{x}-1)
  \\
  &=
  (\cosh{x}-\cos{y})^2,
\end{split}
\end{equation*}
$\cosh{x}\ge1$, and $\cos{y}\le1$, we have $|\cosh{a}-1|=\cosh{x}-\cos{y}$.
\end{proof}
Especially we have $|x|<\arccosh{3/2}=\log\left((3+\sqrt{5})/2\right)=0.9642\dots<1$.
\begin{cor}\label{cor}
If a complex number $a$ satisfies $|2\cosh{a}-2|<1$ and $|\Im{a}|<\pi/3$, then for any real number $u$ with $0<u<1$, we have $|\cosh{ua}-1|<|\cosh{a}-1|$.
\end{cor}
\begin{proof}
From Lemma~\ref{lem:cosh}, $\cosh{x}-\cos{y}<1/2$ with $a:=x+y\sqrt{-1}$.
Since $\cosh{x}$ is increasing (decreasing, respectively) for $x>0$ ($x<0$, respectively) and $\cos{y}$ is decreasing (increasing, respectively) for $0<y<\pi/3$ ($0>y>-\pi/3$, respectively), we have
\begin{multline*}
  |\cosh{ua}-1|
  =
  |\cosh{ux}-\cos{uy}|
  =
  \cosh{ux}-\cos{uy}
  \\
  <
  \cosh{x}-\cos{y}
  =
  |\cosh{x}-\cos{y}|
  =
  |\cosh{a}-1|.
\end{multline*}
\end{proof}
\begin{lem}\label{lem:cosh_ua_va}
For a complex number $a$ with $|2\cosh{a}-2|<1$ and $|\Im{a}|<\frac{\pi}{3}$, and real numbers $u$ and $v$ with $0\le u<v<1$, we have
\begin{equation*}
  |\cosh{a}-\cosh{ua}|\ge|\cosh{a}-\cosh{va}|.
\end{equation*}
Moreover the equality holds only when $a=0$.
\end{lem}
\begin{proof}
It is clear that both hand sides are equal when $a=0$.
So we assume $a\ne0$ and prove the strict inequality.
\par
Put $a:=x+y\sqrt{-1}$ with $(x,y)\ne(0,0)$ and $|y|<\pi/3$.
\par
We will show that
$\varphi(x,y,u):=\left|\cosh{a}-\cosh{ua}\right|^2$
is decreasing with respect to $u$ for $0<u<1$.
Since $\varphi(x,y,u)=\varphi(-x,y,u)=\varphi(x,-y,u)$, we may assume that $x\ge0$ and $\pi/3>y\ge0$.
Since
\begin{equation*}
  |\cosh{a}-\cosh{ua}|^2
  =
  (\cosh{x}\cos{y}-\cosh{ux}\cos{uy})^2
  +
  (\sinh{x}\sin{y}-\sinh{ux}\sin{uy})^2,
\end{equation*}
we have
\begin{equation*}
\begin{split}
  &\frac{\partial\,\varphi(x,y,u)}{\partial\,u}
  \\
  =&
  -2x
  \left\{
    \sinh{x}\cosh{ux}\sin{y}\sin{uy}
    +
    \cosh{x}\sinh{ux}\cos{y}\cos{uy}
    -
    \sinh{ux}\cosh{ux}
  \right\}
  \\
  &-2y
  \left\{
    \sin{uy}\cos{uy}
    +
    \sinh{x}\sinh{ux}\sin{y}\cos{uy}
    -
    \cosh{x}\cosh{ux}\cos{y}\sin{uy}
  \right\}.
\end{split}
\end{equation*}
Put
\begin{align*}
  \varphi_1(x,y,u)
  :=&
  \sinh{x}\cosh{ux}\sin{y}\sin{uy}
  +
  \cosh{x}\sinh{ux}\cos{y}\cos{uy}
  \\
  &-
  \sinh{ux}\cosh{ux}
  \\
  \intertext{and}
  \varphi_2(x,y,u)
  :=&
  \sin{uy}\cos{uy}
  +
  \sinh{x}\sinh{ux}\sin{y}\cos{uy}
  \\
  &-
  \cosh{x}\cosh{ux}\cos{y}\sin{uy}.
\end{align*}
We will show
\begin{enumerate}
\item
  $\varphi_1(x,y,u)>0$ when $x>0$ and $y\ge0$, and
\item
  $\varphi_2(x,y,u)>0$ when $x\ge0$ and $y>0$.
\end{enumerate}
\par
First we will show (1).
Note that if $x>0$, $\varphi_1(x,0,u)=\sinh{ux}(\cosh{x}-\cosh{ux})>0$, and so we will assume that $y>0$.
\par
Since $\varphi_1(0,y,u)=0$, it is sufficient to show that
\begin{equation*}
\begin{split}
  &\frac{\partial\,\varphi_1(x,y,u)}{\partial\,x}
  \\
  =&
  u(\cosh{x}\cosh{ux}\cos{y}\cos{uy}+\sinh{x}\sinh{ux}\sin{y}\sin{uy}
  -\sinh^2{ux}-\cosh^2{ux})
  \\
  &+
  \cosh{x}\cosh{ux}\sin{y}\sin{uy}+\sinh{x}\sinh{ux}\cos{y}\cos{uy}
\end{split}
\end{equation*}
is positive when $x>0$, $\pi/3>y>0$, and $1>u>0$.
Note that
\begin{equation*}
  \frac{\partial\,\varphi_1(x,y,u)}{\partial\,x}\bigr|_{x=0}
  =
  u(\cos{y}\cos{uy}-1)+\sin{y}\sin{uy}
\end{equation*}
is positive since its partial derivative with respect to $y$ is
$\left(1-u^2\right)\cos{y}\sin{uy}$, which is positive.
Moreover we have
\begin{equation*}
\begin{split}
  &\frac{\partial^2\,\varphi_1(x,y,u)}{\partial\,x^2}
  \\
  =&
  \cosh{ux}
  \left[
    \sinh{x}
    \left\{
      2u\cos{y}\cos{uy}+\left(1+u^2\right)\sin{y}\sin{uy}
    \right\}
    -2u^2\sinh{ux}
  \right]
  \\
  &+
  \sinh{ux}
  \left[
    \cosh{x}
    \left\{
      2u\sin{y}\sin{uy}+\left(1+u^2\right)\cos{y}\cos{uy}
    \right\}
    -2u^2\cosh{ux}
  \right]
  \\
  >&
  \cosh{ux}
  \left[
    \sinh{x}
    \left\{
      2u\cos{y}\cos{uy}+\left(1+u^2\right)\sin{y}\sin{uy}
    \right\}
    -2u^2\sinh{x}
  \right]
  \\
  &+
  \sinh{ux}
  \left[
    \cosh{x}
    \left\{
      2u\sin{y}\sin{uy}+\left(1+u^2\right)\cos{y}\cos{uy}
    \right\}
    -2u^2\cosh{x}
  \right]
  \\
  =&
  \cosh{ux}
  \sinh{x}
  \left\{
    2u\cos{y}\cos{uy}+\left(1+u^2\right)\sin{y}\sin{uy}-2u^2
  \right\}
  \\
  &+
  \sinh{ux}
  \cosh{x}
  \left\{
    2u\sin{y}\sin{uy}+\left(1+u^2\right)\cos{y}\cos{uy}-2u^2
  \right\}
  \\
  =&
  \sinh{x}
  \cosh{ux}
  \left\{
    2u\cos{y}\cos{uy}+(1-u)^2\sin{y}\sin{uy}+2u\sin{y}\sin{uy}-2u^2
  \right\}
  \\
  &+
  \sinh{ux}
  \cosh{x}
  \left\{
    2u\sin{y}\sin{uy}+(1-u)^2\cos{y}\cos{uy}+2u\cos{y}\cos{uy}-2u^2
  \right\}
  \\
  >&
  (\sinh{x}\cosh{ux}+\sinh{ux}\cosh{x})
  \left(2u\sin{y}\sin{uy}+2u\cos{y}\cos{uy}-2u^2\right)
  \\
  =&
  2u(\sinh{x}\cosh{ux}+\sinh{ux}\cosh{x})
  (\cos{(1-u)y}-u)
  \\
  >&
  2u(\sinh{x}\cosh{ux}+\sinh{ux}\cosh{x})
  \left\{-\frac{3}{2\pi}(1-u)y+1-u\right\}
  \\
  =&
  2u(1-u)(\sinh{x}\cosh{ux}+\sinh{ux}\cosh{x})
  \left(1-\frac{3y}{2\pi}\right)
  \\
  >&
  0,
\end{split}
\end{equation*}
since $0<u<1$ and $\cos{z}>-\dfrac{3}{2\pi}z+1$ for $0<z<\pi/3$.
Therefore $\partial\,\varphi_1(x,y,u)/\partial\,x$ is also positive.
\par
Next we will show (2).
Note that if $\pi/3>y>0$, $\varphi_2(0,y,u)=\sin{uy}(\cos{uy}-\cos{y})>0$, and so we will assume that $x>0$.
Since $\varphi_2(x,y,0)=0$, it is sufficient to show that
\begin{equation*}
\begin{split}
  \frac{\partial\,\varphi_2(x,y,u)}{\partial\,x}
  =&
  \cosh{x}\sinh{ux}
  \left(
    \sin{y}\cos{uy}
    -
    u\sin{uy}\cos{y}
  \right)
  \\
  &+
  \sinh{x}\cosh{ux}
  \left(
    u\sin{y}\cos{uy}
    -
    \sin{uy}\cos{y}
  \right)
\end{split}
\end{equation*}
is positive when $x>0$, $\pi/3>y>0$, and $1>u>0$.
The first term is clearly positive and so we will show that
$u\sin{y}\cos{uy}-\sin{uy}\cos{y}$ is positive.
But this can be easily verified since it is $0$ when $y=0$ and its derivative with respect to $y$ is $\left(1-u^2\right)\sin{y}\sin{uy}$, which is positive.
\end{proof}
\begin{lem}\label{lem:cosh_ua}
There exists a positive number $\varepsilon$ such that for a complex number $a\ne0$ with $|\Im{a}|<\frac{\pi}{3}$ and $|\Re{a}|<\pi$, and a real number $u$ with $0<u<\varepsilon$, we have
\begin{equation*}
  \left|
    \frac{\cosh{a}-\cosh{ua}}{\cosh{a}-1}
  \right|
  >
  1-u.
\end{equation*}
\end{lem}
\begin{proof}
We will show that
\begin{equation*}
  |\cosh{a}-\cosh{ua}|^2-(1-u)^2|\cosh{a}-1|^2>0
\end{equation*}
if $0<u<\varepsilon$.
Putting $a:=x+y\sqrt{-1}$ with $|x|<\pi$ and $|y|<\frac{\pi}{3}$, the left hand side equals
\begin{equation*}
\begin{split}
  &(\cosh{x}\cos{y}-\cosh{ux}\cos{uy})^2
  +
  (\sinh{x}\sin{y}-\sinh{ux}\sin{uy})^2
  \\
  &-
  (1-u)^2
  \{(\cosh{x}\cos{y}-1)^2+\sinh^2{x}\sin^2{y}\}
  \\
  =&
  \left\{(2-u)\sinh{x}\sin{y}-\sinh{ux}\sin{uy}\right\}
  \\
  &\quad\times
  \left(u\sinh{x}\sin{y}-\sinh{ux}\sin{uy}\right)
  \\
  &+
  \left\{u-1+(2-u)\cosh{x}\cos{y}-\cosh{ux}\cos{uy}\right\}
  \\
  &\quad\times
  \left(1-u+u\cosh{x}\cos{y}-\cosh{ux}\cos{uy}\right).
\end{split}
\end{equation*}
Since it remains the same if we alter the signs of $x$ or $y$, we may assume that $\pi>x\ge0$ and $\pi/3>y\ge0$ ($(x,y)\ne(0,0)$).
Put
\begin{align*}
  \alpha_1(x,y,u)
  &:=
  (2-u)\sinh{x}\sin{y}-\sinh{ux}\sin{uy},
  \\
  \alpha_2(x,y,u)
  &:=
  u\sinh{x}\sin{y}-\sinh{ux}\sin{uy},
  \\
  \beta_1(x,y,u)
  &:=
  u-1+(2-u)\cosh{x}\cos{y}-\cosh{ux}\cos{uy},
  \\
  \beta_2(x,y,u)
  &:=
  1-u+u\cosh{x}\cos{y}-\cosh{ux}\cos{uy}.
\end{align*}
We will show that $\alpha_1(x,y,u)$, $\alpha_2(x,y,u)$, $\beta_1(x,y,u)$, and $\beta_2(x,y,u)$ are all positive.
\par
Since $0<u<1$, $\sinh{x}$ is increasing for any $x$, and $\sin{y}$ is increasing when $0\le y<\pi/3$, we have
\begin{equation*}
  \alpha_1(x,y,u)
  >
  (2-u)\sinh{ux}\sin{uy}-\sinh{ux}\sin{uy}
  =
  (1-u)\sinh{ux}\sin{uy}
  >0.
\end{equation*}
\par
By the Taylor expansion of $\alpha_2(x,y,u)$ around $x=0$, we have
\begin{equation*}
  \alpha_2(x,y,u)
  =
  \sum_{n=0}^{\infty}\frac{1}{(2n+1)!}u(\sin{y}-u^{2n}\sin{uy})x^{2n+1}.
\end{equation*}
Since $\sin{y}$ is increasing for $0<y<\pi/3$ and $0<u<1$, we have $\sin{y}-u^{2n}\sin{uy}>0$.
Therefore $\alpha_2(x,y,u)>0$.
\par
The Taylor expansions of $\beta_1(x,y,u)$ and $\beta_2(x,y,u)$ around $u=0$ give
\begin{align*}
  \beta_1(x,y,u)
  &=
  2(\cosh{x}\cos{y}-1)
  -
  (\cosh{x}\cos{y}-1)u
  +
  \frac{1}{2}(y^2-x^2)u^2
  \\
  &\quad
  +
  \frac{1}{24}(-x^4+6x^2y^2-y^4)u^4
  +
  O\left(u^6\right)
  \\
  \beta_2(x,y,u)
  &=
  (\cosh{x}\cos{y}-1)u
  +
  \frac{1}{2}(y^2-x^2)u^2
  \\
  &\quad
  +
  \frac{1}{24}
  (-x^4+6x^2y^2-y^4)u^4
  +
  O\left(u^6\right)
\end{align*}
and so
\begin{multline*}
  \beta(x,y,u)
  :=
  \beta_1(x,y,u)\beta_2(x,y,u)
  \\
  =
  2(\cosh{x}\cos{y}-1)^2u
  -
  (\cosh{x}\cos{y}-1)(\cosh{x}\cos{y}-1+x^2-y^2)u^2
  \\
  +
  \frac{1}{12}
  \left(
    4(x^4-3x^2y^2+y^4)-(x^4-6x^2y^2+y^4)\cosh{x}\cos{y}
  \right)
  u^4
  +
  O\left(u^6\right).
\end{multline*}
Therefore $\beta(x,y,u)$ is positive for small $u$ if $\cosh{x}\cos{y}\ne1$.
\par
When $\cosh{x}\cos{y}=1$ we have
\begin{equation}\label{eq:beta}
  \beta(x,y,u)
  =
  \frac{1}{4}(x^2-y^2)^2u^4+O\left(u^6\right).
\end{equation}
Since
\begin{align*}
  \cosh{x}\cos{x}\big|_{x=0}
  &=1,
  \\
  \left.\frac{d\,\cosh{x}\cos{x}}{d\,x}\right|_{x=0}
  &=
  (\sinh{x}\cos{x}-\cosh{x}\sin{x})\big|_{x=0}
  =0,
  \intertext{and}
  \\
  \frac{d^2\,\cosh{x}\cos{x}}{d\,x^2}
  &=
  -2\sinh{x}\sin{x}<0\quad\text{if $0<x<\pi$},
\end{align*}
$\cosh{x}\cos{x}<1$ for $0<x<\pi$, which means that $x\ne{y}$ when $\cosh{x}\cos{y}=1$.
So $\beta(x,y,u)>0$ for small $u$ since the coefficient of $u^4$ in \eqref{eq:beta} is positive.
\par
Thus we have concluded that $\beta(x,y,u)>0$ for small $u$.
\end{proof}
\begin{lem}\label{lem:ExpIntegralE}
For a positive integer $m$ and a positive real number $a$, we have
\begin{equation*}
  \int_{1}^{\infty}e^{-at}t^{m}dt
  =
  \frac{e^{-a}}{a}\sum_{k=0}^{m}\frac{m!}{a^{k}(m-k)!}
\end{equation*}
\end{lem}
\begin{proof}
Integration by parts gives
\begin{equation*}
\begin{split}
  &
  \int_{1}^{\infty}e^{-at}t^{m}dt
  \\
  =&
  \left[-\frac{1}{a}e^{-at}t^{m}\right]_{1}^{\infty}
  +
  \frac{m}{a}\int_{1}^{\infty}e^{-at}t^{m-1}dt
  \\
  =&
  \frac{1}{a}e^{-a}
  +
  \frac{m}{a}\int_{1}^{\infty}e^{-at}t^{m-1}dt
  \\
  =&
  \frac{1}{a}e^{-a}
  +
  \frac{m}{a^2}e^{-a}
  +
  \frac{m(m-1)}{a^2}\int_{1}^{\infty}e^{-at}t^{m-2}dt
  \\
  =&
  \frac{1}{a}e^{-a}
  +
  \frac{m}{a^2}e^{-a}
  +\dots+
  \frac{m(m-1)\times\dots\times(m-k+1)}{a^{k+1}}e^{-a}
  +\cdots
  \\
  &+
  \frac{m(m-1)\times\dots\times2}{a^{m}}e^{-a}
  +
  \frac{m!}{a^{m}}\int_{1}^{\infty}e^{-at}dt
  \\
  =&
  \frac{1}{a}e^{-a}
  +
  \frac{m}{a^2}e^{-a}
  +\dots+
  \frac{m(m-1)\times\dots\times(m-k+1)}{a^{k+1}}e^{-a}
  +\cdots
  \\
  &+
  \frac{m(m-1)\times\dots\times2}{a^{m}}e^{-a}
  +
  \frac{m!}{a^{m+1}}e^{-a}
\end{split}
\end{equation*}
and the proof is complete.
\end{proof}
\begin{lem}\label{lem:cosh_x}
For a complex number $a$ with $|2\cosh{a}-2|<1$, $|\Im{a}|<\pi/3$, and $a\ne0$,
there exists a positive number $\varepsilon'>0$ such that if $0<x<\varepsilon'$ and $0<u<\varepsilon'$, then
\begin{equation*}
  1>
  \left|
    \frac{\cosh{a(1-x)}-\cosh{ua}}
         {\cosh{a}-\cosh{ua}}
  \right|
  >
  1-
  \left|
    \frac{a\sinh{a}}{\cosh{a}-1}
  \right|x.
\end{equation*}
\end{lem}
\begin{proof}
Using the Taylor expansion with respect to $x$ and $u$ around $x=u=0$, we have
\begin{equation*}
  \frac{\cosh{a(1-x)}-\cosh{ua}}
       {\cosh{a}-\cosh{ua}}
  =
  1
  -
  \frac{a\sinh{a}}{\cosh{a}-1}x
  +
  \frac{a^2\cosh{a}}{2(\cosh{a}-1)}x^2
  +
  R_a(u,x),
\end{equation*}
where $R_a(u,x)$ is the terms with total degrees of $u$ and $x$ are greater than two.
\par
From Lemma~\ref{lem:Re_a}, we have
\begin{equation*}
  \left|
    \frac{\cosh{a(1-x)}-\cosh{ua}}
         {\cosh{a}-\cosh{ua}}
  \right|
  <
  1
\end{equation*}
if $x$ and $u$ are sufficiently small.
\par
Moreover we have
\begin{equation*}
\begin{split}
  &
  \left|
    \frac{\cosh{a(1-x)}-\cosh{ua}}
       {\cosh{a}-\cosh{ua}}
  \right|
  +
  \left|
    \frac{a\sinh{a}}{\cosh{a}-1}
  \right|
  x
  \\
  &\quad\ge
  \left|
    \frac{\cosh{a(1-x)}-\cosh{ua}}
       {\cosh{a}-\cosh(ua)}
    +
    \frac{a\sinh{a}}{\cosh{a}-1}x
  \right|
  \\
  &\quad=
  \left|
    1
    +
    \frac{a^2\cosh{a}}{2(\cosh{a}-1)}x^2
  \right|
  +
  R_a(u,x).
\end{split}
\end{equation*}
From Lemma~\ref{lem:Re_a^2}, we have $\Re{\frac{a^2\cosh{a}}{2(\cosh{a}-1)}}>0$ if $|2\cosh{a}-2|<1$, and the other inequality follows.
\end{proof}
\begin{lem}\label{lem:Re_a}
For a complex number $a\ne0$ with $|\Im{a}|<\pi$, we have
\begin{equation*}
  \Re\frac{a\sinh{a}}{\cosh{a}-1}>0.
\end{equation*}
\end{lem}
\begin{proof}
We put $a:=x+\sqrt{-1}y$ with $|y|<\pi$.
Then we have
\begin{equation*}
  \Re\frac{a\sinh{a}}{\cosh{a}-1}
  =
  \frac{(\cosh{x}-\cos{y})(x\sinh{x}+y\sin{y})}
       {(\cos{y}\cosh{x}-1)^2+\sin^2{y}\sinh^2{x}}
  >0
\end{equation*}
if $|y|<\pi$ and $(x,y)\ne(0,0)$.
\end{proof}
\begin{lem}\label{lem:Re_a^2}
For a complex number $a\ne0$ with $|2\cosh{a}-2|<1$ and $|\Im{a}|<\pi/3$, we have
\begin{equation*}
  \Re{\frac{a^2\cosh{a}}{\cosh{a}-1}}>0.
\end{equation*}
\end{lem}
\begin{proof}
Putting $a:=x+\sqrt{-1}y$ with $x,y\in\R$, we have
\begin{equation*}
  \Re{\frac{a^2\cosh{a}}{\cosh{a}-1}}
  =
  \frac{f(x,y)}
  {(\cos{y}\cosh{x}-1)^2+\sin^2{y}\sinh^2{x}}
\end{equation*}
with
\begin{equation*}
  f(x,y)
  :=
  (x^2-y^2)(\cos^2{y}+\sinh^2{x}-\cos{y}\cosh{x})
  +
  2xy\sin{y}\sinh{x}.
\end{equation*}
\par
We will show that $f(x,y)>0$.
We may assume $x\ge0$, $y\ge0$ ($(x,y)\ne(0,0)$) as before.
Since
\begin{equation*}
  f(x,0)
  =
  x^2(1+\sinh^2{x}-\cosh{x})
  =
  x^2(\cosh^2{x}-\cosh{x})
  >0
\end{equation*}
when $x>0$, we will assume that $y>0$.
\par
Since $f(x,y)$ is analytic, it is sufficient to prove that every $n$th derivative of $f$ at $x=0$ is positive when $n$ is even and zero when $n$ is odd.
\par
Since
\begin{align*}
  \frac{\partial^k\,(x^2-y^2)}{\partial\,x^k}\Biggr|_{x=0}
  &=
  \begin{cases}
  -y^2
  &\quad\text{if $k=0$,}
  \\
  2
  &\quad\text{if $k=2$,}
  \\
  0
  &\quad\text{otherwise,}
  \end{cases}
  \\
  \frac{\partial^k\,(\sinh^2{x}-\cos{y}\cosh{x})}
       {\partial\,x^k}\Biggr|_{x=0}
  &=
  \begin{cases}
  -\cos{y}
  &\quad\text{if $k=0$,}
  \\
  2^{k-1}-\cos{y}
  &\quad\text{if $k$ is even and positive,}
  \\
  0
  &\quad\text{otherwise,}
  \end{cases}
  \\
  \intertext{and}
  \frac{\partial^k\,(x\sinh{x})}{\partial\,x^k}
  &=
  \begin{cases}
  k
  &\quad\text{if $k$ is even and positive,}
  \\
  0
  &\quad\text{otherwise,}
  \end{cases}
\end{align*}
we have
\begin{equation*}
\begin{split}
  &\frac{\partial^n\,f(x,y)}{\partial\,x^n}\Biggr|_{x=0}
  \\
  =&
  \sum_{k=0}^{n}
  \binom{n}{k}
  \frac{\partial^{n-k}\,(x^2-y^2)}{\partial\,x^{n-k}}\Biggr|_{x=0}
  \times
  \frac{\partial^k(\cos^2{y}+\sinh^2{x}-\cos{y}\cosh{x})}
       {\partial\,x^k}\Biggr|_{x=0}
  \\
  &+
  2y\sin{y}\frac{\partial^n\,(x\sinh{x})}{\partial\,x^n}\Biggr|_{x=0}
  \\
  =&
  \begin{cases}
  -y^2(\cos^2{y}-\cos{y})
  &\quad\text{if $n=0$,}
  \\
  -y^2(2-\cos{y})+2(\cos^2{y}-\cos{y})+4y\sin{y}
  &\quad\text{if $n=2$,}
  \\
  -y^2(2^{n-1}-\cos{y})+n(n-1)(2^{n-3}-\cos{y})+2ny\sin{y}
  &\quad\text{if $n>3$ and even,}
  \\
  0
  &\quad\text{otherwise.}
  \end{cases}
\end{split}
\end{equation*}
It is clear that $-y^2(\cos^2{y}-\cos{y})>0$ since $y<\pi/3$.
If $n$ if even and $n\ge4$, then since $y<\pi/3$ we have
\begin{equation*}
\begin{split}
  &\frac{\partial^n\,f(x,y)}{\partial\,x^n}\Biggr|_{x=0}
  \\
  =&
  2^{n-4}
  \left\{
    n(n-1)-8y^2
  \right\}
  +
  n(n-1)\left(2^{n-4}-\cos{y}\right)
  +
  y^2\cos{y}
  +
  2ny\sin{y}
  \\
  >&
  2^{n-4}
  \left\{
    12-8\left(\frac{\pi}{3}\right)^2
  \right\}
  +
  n(n-1)(1-\cos{y})
  >0.
\end{split}
\end{equation*}
To show that
\begin{equation*}
  \frac{\partial^2\,f(x,y)}{\partial\,x^2}\Biggr|_{x=0}
  =
  -2y^2+2(\cos^2{y}-\cos{y})+4y\sin{y}+y^2\cos{y}
\end{equation*}
is positive, we will consider the function
$g(y):=-2y^2+2(\cos^2{y}-\cos{y})+4y\sin{y}$.
Since
\begin{equation*}
  \frac{d\,g(y)}{d\,y}
  =
  2(3\sin{y}-2y)+4\cos{y}(y-\sin{y})
\end{equation*}
is easily verified to be positive, we have $g(y)>0$.
So $\partial^2\,f(x,y)/\partial\,x^2\bigr|_{x=0}$ is also positive.
\end{proof}
\bibliography{mrabbrev,hitoshi}
\bibliographystyle{hamsplain}
\end{document}